\documentclass[12pt]{article}       
\usepackage{amsfonts,amsmath,amssymb}
\newtheorem{Def}{Definition}
\newtheorem{Th}{Theorem}

\newtheorem{lem}[Th]{Lemma}

\newtheorem{Conjecture}{Conjecture}
\newcommand{\D}{{\Delta}}

\newenvironment{Pf}
{\noindent\textbf{Proof.}\ \ }{\hfill $\Box$ \medskip}

\title{The strong equitable vertex 2-arboricity of complete bipartite and tripartite graphs}

\author {Keaitsuda Maneeruk Nakprasit\\ {\small\em Department of Mathematics, Faculty of Science, Khon Kaen University, 40002, Thailand }\\  
{\small\em E-mail address: kmaneeruk@hotmail.com} 
\and Kittikorn Nakprasit \footnote{Corresponding Author} \\ 
{\small\em Department of Mathematics, Faculty of Science, Khon Kaen University, 40002, Thailand }\\
{\small\em E-mail address: kitnak@hotmail.com}}
\date{}

\begin{document}

\maketitle

\begin{center}{\bf Abstract}\end{center}
\indent\indent

A $(q,r)$\emph{-tree-coloring} of a graph $G$ is a $q$-coloring of vertices of $G$   
such that the subgraph induced by each color 
class is a forest of maximum degree at most $r.$ 
An \emph{equitable $(q, r)$-tree-coloring} of a graph $G$ is a $(q,r)$-tree-coloring  
such that the sizes of any two color classes differ by at most one. 
Let the \emph{strong equitable vertex $r$-arboricity} 
be the minimum $p$ such that $G$ has an equitable
$(q, r)$-tree-coloring for every $q\geq p.$ 


In this paper, we find the exact value for each $va^\equiv_2(K_{m,n})$ 
and $va^\equiv_2(K_{l,m,n}).$

\section{Introduction} 
Throughout this paper, all graphs are finite, undirected, and 
simple. We use  $V(G)$ and $E(G),$ respectively, 
to denote the vertex set and edge set of a graph $G.$ 
We use  $V(G)$ and $E(G),$ respectively, 
to denote the vertex set and edge set of a graph $G.$ 
For a complete bipartite graph $K_{m,n}$ where $m \leq n,$ 
we let $X=\{x_1,\ldots,x_m\}$ and 
$Y=\{y_1,\ldots,y_n\}$ to be the partite sets of $K_{m,n}.$  
For a complete tripartite graph $K_{l,m,n}$ where $l \leq m \leq n,$ 
we have $X=\{x_1,\ldots,x_l\},$   
$Y=\{y_1,\ldots,y_m\},$ and  $Z=\{z_1,\ldots,z_n\},$ to be 
the partite sets of $K_{l,m,n}.$

An \emph{equitable $k$-coloring} of a graph is a proper vertex
$k$-coloring such that the sizes of every two color classes differ
by at most $1.$ 

It is known \cite{GareyJohnson} that determining if a planar graph
with maximum degree $4$ is $3$-colorable is NP-complete. For a given
$n$-vertex planar graph $G$ with maximum degree $4,$ let $G'$ be the
graph obtained from $G$ by adding $2n$ isolated vertices. Then $G$
has $3$-coloring if and only if $G'$ has an equitable $3$-coloring.
Thus, finding the minimum number of colors needed to color a graph
equitably even for a planar graph is an NP-complete problem.

Hajnal and Szemer\'edi~\cite{HS} settled a conjecture of Erd\H os by
proving that  every graph $G$ with maximum degree at most $\Delta$
has an equitable $k$-coloring for every $k\geq 1+\Delta.$ 
This result is now known as Hajnal and Szemer\'edi Theorem. 
Later, Kierstead and
Kostochka~\cite{KK08} gave a simpler proof of Hajnal and Szemer\'edi
Theorem. The bound of the
Hajnal-Sz{e}mer\' edi theorem is sharp, but it can be improved for
some important classes of graphs. In fact, Chen, Lih, and
Wu~\cite{CLW94} put forth the following conjecture.

\begin{Conjecture} \label{ConjLW}
Every connected graph $G$ with maximum degree $\Delta\geq 2$ has an
equitable coloring with $\Delta$ colors, except when $G$ is a
complete graph or an odd cycle or $\Delta$ is odd and
$G=K_{\Delta,\Delta}.$
\end{Conjecture}

Lih and Wu~\cite{LW} proved the conjecture for bipartite graphs.
Meyer \cite{M} proved that every forest with maximum degree $\Delta$
has an equitable $k$-coloring for each $k \geq 1+\lceil
\Delta/2\rceil $ colors. This result implies the conjecture holds
for forests.  Yap and Zhang~\cite{YZ1} proved that the
conjecture holds for outerplanar graphs. Later Kostochka~\cite{Ko}
improved the result  by proving that every
outerplanar graph with maximum degree $\Delta$ has an equitable
$k$-coloring for each $k \geq 1+\lceil \Delta/2\rceil.$ 

In~\cite{ZY98}, Zhang and Yap essentially proved the conjecture
holds for planar graphs with maximum degree at least $13.$ Later
Nakprasit~\cite{Nak12} extended the result to all planar graphs with
maximum degree at least  $9.$ 
Some related results are about planar graphs without some restricted cycles 
~\cite{LiBu09, NN12, ZhuBu08}. 

Moreover, the conjecture has been confirmed for other classes of graphs, 
such as graphs with degree at most 3~\cite{CLW94, CY12}  
and series-parallel graphs \cite{ZW11}. 

In contrast with ordinary coloring, a graph may have an equitable
$k$-coloring  but has no equitable $(k+1)$-coloring. 
For example, $K_{7,7}$ has an equitable $k$-coloring
for $k=2,4,6$ and $k \ge 8$, but has no equitable $k$-coloring for
$k=3,5$ and $7$. This leads to the definition of the 
\emph{equitable chromatic threshold} which is 
is the minimum $p$ such that $G$ has an equitable
$q$-coloring for every $q\geq p,$

In \cite{Fan11}, Fan, Kierstead, Liu, Molla, Wu, and Zhang 
considered an equitable relaxed colorings. 
They proved that every graph with maximum degree $\D$ has an equitable $\D$-coloring.  
such that each color class induces a forest with maximum degree at most one. 

On the basis of the aforementioned research,  Wu, Zhang, and Li \cite{WZL13} introduced 
a $(q,r)$\emph{-tree-coloring} of a graph $G$ which is a $q$-coloring of vertices of $G$   
such that the subgraph induced by each color class is a forest of maximum degree at most $r.$ 
A $(q,\infty)$\emph{-tree-coloring} of a graph $G$ is a $q$-coloring of $G$ such that  
the subgraph induced by each color class is a forest. 
An \emph{equitable $(q, r)$-tree-coloring} of a graph $G$ is a $(q,r)$-tree-coloring  such that the
sizes of any two color classes differ by at most one. 
Thus, the result of Fan, Kierstead, Liu, Molla, Wu, and Zhang can be restated that  
every graph with maximum degree $\D$ has an equitable $(\D,1)$-tree-coloring.  

Let the \emph{strong equitable vertex $k$-arboricity}, denoted by $va^\equiv_r (G),$ 
be the minimum $p$ such that $G$ has an equitable
$(q, r)$-tree-coloring for every $q\geq p.$ 
Wu, Zhang, and Li \cite{WZL13} 
 proved  that $va^\equiv_\infty(G) \leq 3$ for each planar graph $G$ with girth at least 5 
and $va^\equiv_\infty(G) \leq 3$ for each  planar graph $G$ with girth 
at least 6 and for each outerplanar graph. 
Moreover, they gave a sharp upper bound for $va^\equiv_1(K_{n,n})$ in general case.  
They commented that finding the strong equitable 1-arboricity 
for every $K_{n,n}$ seems not to be an easy task.


In this paper, we find the exact value for each $va^\equiv_2(K_{m,n})$ 
and $va^\equiv_2(K_{l,m,n}).$   

\section{Useful Lemmas} 
We introduce the notion of $p(q: n_1,\ldots, n_k)$ which can be computed 
in  linear-time.   

\begin{Def}
Assume that $G=K_{n_1,\ldots,n_k}$ has an  equitable $q$-coloring,  
and $d$ is the minimum value greater than $\lfloor (n_1+\cdots +n_k)/q \rfloor$ 
such that (i) there are distinct $i$ and $j$ in which $n_i$ and $n_j$ 
are not divisible by $d,$ or 
(ii) there is $n_j$ with  $n_j/\lfloor n_j/d\rfloor  > d+1.$ 
Define $p(q: n_1,n_2,\ldots, n_k)= \lceil n_1/d \rceil+\cdots+\lceil n_k/d \rceil.$  
\end{Def}

\begin{Th}\label{T1} \cite{NN1ab}
Assume that $G=K_{n_1,\ldots,n_k}$ has an  equitable $q$-coloring.   
Then $p(q: n_1,\ldots, n_k)$ is the minimum $p$ such that 
$G$ is equitable $r$-colorable for each $r$ satisfying $p \leq r \leq q.$ 
\end{Th} 

\begin{lem} \label{proper} 
Let $G$ be a complete multipartite graph with $n$ vertices. 
If the size of a color class from a $(q,2)$-tree coloring of $G$ 
is at least 4, 
then the color class is independent. 
Consequently,  each equitable $(q,2)$-tree coloring of $G$ 
such that $n/q \geq 4$ is a proper equitable coloring.  
\end{lem} 
\begin{Pf}
Suppose to the contrary that  a color class $C$ of of size $k\geq 4$ 
is not an independent set. 
Then $C$ induces $K_{1,k-1}$ or a graph with a cycle, a contradiction. 
The remaining of the Lemma follows immediately. 
\end{Pf}

\begin{lem}\label{L3} 
Let $G=K_{n_1,\ldots,n_k}$ and $N=n_1+\cdots+n_k.$ 
Assume $G$ has an equitable $q$-coloring where $N/(q-1) \geq 4$ and 
$G$ has an equitable $(r,2)$-tree-coloring for each $r \geq q.$ 
Then $va^\equiv_2(G) = p(q: n_1,\ldots,n_k).$ 
\end{lem}
\begin{Pf} 
Let $p =p(q: n_1,\ldots,n_k).$ 
From the definition of $p$ and the condition of $q,$ 
the graph $G$ has an equitable $(r,1)$-tree-coloring for each $r \geq p.$ 
To complete the proof,  
it suffices to show that $G$ has no equitable $(p-1,2)$-tree-coloring. 
Suppose to the contrary that $G$ has an equitable $(p-1,2)$-tree-coloring. 
Since $p-1 \leq  q-1,$ 
each color class has size at least $n/(p-1) \geq n/(q-1) \geq 4.$ 
Lemma \ref{proper} yields that $G$ has a proper equitable $(p-1)$-coloring. 
But this 

\end{Pf}

Let $G =K_{m,n}$ or $K_{l,m,n}.$ 
We introduce an algorithm to construct a $(q,2)$-tree-coloring of $G.$ 
The first key idea is that we arrange vertices of $G$ in a way that vertices 
in a same partite set are consecutively ordered. 
Then we partition $V(G)$ in a way that each partitioned set (color class) 
contains $k$ or $k+1$ consecutive vertices from the  arrangement. 
By this method, there are at most one non-independent color class in $K_{m,n},$ 
and at most two non-independent color classes in $K_{l,m,n}.$ 

The second key idea is that we want each non-independent color class 
to have size at most 3. 
The final key idea is that we want elements 
in each non-independent color class comes from exactly two partite sets. 
To achieve this objective for any $K_{l,m,n}$ except  $K_{1,1,1},$ 
we have elements in $Z$ (a largest partite set with size at most 2) 
in the middle of the arrangement. 

A coloring satisfying all of these three key ideas has 
each non-independent color class induces a tree of maximum degree at most 2. 
If the sizes of any two color classes differ by at most one, 
then we have an equitable $(q,2)$-tree-coloring. 
Now we show a desired algorithm to obtain an equitable $(q,2)$-tree-coloring as follows.  
 

\begin{Def} (Algorithm A) Let $G =K_{m,n}$ or $K_{l,m,n}$ and $k\leq 3.$  
If $G=K_{m,n},$ then we let $(v_1,\ldots,v_{m+n})=(x_1,\ldots,x_m,$ $y_1,\ldots,y_n),$ 
otherwise $(v_1,\ldots,v_{l+m+n})$ $=$ $(x_1,\ldots,x_l,z_1,$ 
$\ldots,z_n,$ $y_1,\ldots,y_m).$ 
Consider $|V(G)|=r_0k+s_0(k+1)$ where $r_0$ is a positive integer 
and $k$ and $s_0$ are  nonnegative integers.\\
(1) Set $i=1$ and $j=0.$  \\
(2) 
If $s_j\geq 1$ and (i) $k \leq 2$ or (ii) $v_i, v_{i+k}$ are in a same partite set, 
then let $A_{j+1}= \{v_i,\ldots, v_{i+k}\},$ $r_{j+1}=r_j,$ $s_{j+1}=s_j-1$ 
and set $i\leftarrow i+k+1.$ 
Otherwise, let $A_{j+1}= \{v_i,\ldots, v_{i+k-1}\},$ $r_{j+1}=r_j-1,$ $s_{j+1}=s_j$ 
and set $i\leftarrow i+k.$ \\ 
(3) Set $j \leftarrow j+1.$ 
If $r_j+s_j =\geq 1,$ then go to step (2), otherwise stop.    
\end{Def}

Note that we use Algorithm A for $k \leq 3.$ 
If $k=3,$ then we require $r_0 \geq 1$ for $K_{m,n},$ 
and $r_0 \geq 2$ for $K_{l,m,n}.$

\section{$va^\equiv_2(K_{m,n})$} 

\begin{lem}\label{bi0} 
Let $m+n =4b+c$ where $b$ is a nonnegative integer and $0 \leq c \leq 3.$  
Then $K_{m,n}$ has an equitable $(t,2)$-tree-coloring for each $t\geq b+1.$ 
\end{lem} 
\begin{Pf} 
Let $m+n=4b+c =rk+s(k+1)$ where $r$ is a positive integer,  
$s$ and $k$ are nonnegative integers. 
First consider the case $r+s=b+1.$ 
Then $k \leq 3.$ 
If $k\leq 2,$ then each color class from Algorithm A is 
an independent set or induces $K_1,K_2,$ or $K_{1,2}.$ 
Thus, we obtain an equitable $(r+s=b+1,2)$-tree-coloring. 

Now, we assume $k =3.$ 
Consequently, $c=0,1,2,$ or $3.$ 
If $c=0,$ then $r=4, s=b-3,$ and $b \geq 3.$ 
If $c=1,$ then $r=3, s=b-2,$ and $b \geq 2.$  
If $c=2,$ then $r=2, s=b-1,$ and $b \geq 1.$ 
If $c=3,$ then $r=1$ and  $s=b.$ 

We show that Algorithm A yields an equitable $(b+1,2)$-tree coloring. 
By step (2) of Algorithm A, a non-independent color class (if exists) 
contains three elements from two partite sets. 
Then each color class is independent set or induces $K_{1,2}.$ 
Thus we obtain an equitable $(b+1,2)$-tree-coloring. 

Finally, consider the case that $r+s \geq b+2.$ 
Again we have  (i) $k\leq 2$ or (ii) $k=3$ and $r\geq 2.$   
Similar to the case of $r+s=b+1,$ 
we can use Algorithm A  to 
obtain  an equitable $(r+s,2)$-tree-coloring. 
This completes the proof. 
\end{Pf}

\begin{lem}\label{bi1} 
Let $m+n =4b+c$ where $b$ is a nonnegative  integer and $0\leq c \leq 3.$ 
If $K_{m,n}$ has an equitable $(b,2)$-tree-coloring, 
then  $va^\equiv_2(K_{m,n})=p(b: m,n),$ 
otherwise  $va^\equiv_2(K_{m,n})=b+1.$\\ 
\end{lem}
\begin{Pf} 
From Lemma \ref{bi0}, 
$K_{m,n}$ has an equitable $(t,2)$-tree-coloring for each $t\geq b+1.$ 

If $K_{m,n}$ has no equitable $(b,2)$-tree-coloring, then  
 $va^\equiv_2(K_{m,n})=b+1$ 
by definition of $va^\equiv_2(K_{m,n}).$

Assume $K_{m,n}$ has an equitable $(b,2)$-tree-coloring. 
Then each color class has size at least 4. 
By Lemma \ref{proper}, such equitable  $(b,2)$-tree-coloring 
is a proper equitable $b$-coloring. 
Thus $va^\equiv_2(K_{m,n})=p(b: m,n)$ by Lemma \ref{L3}. 
\end{Pf}

\begin{Th}\label{biM} 
$va^\equiv_2(K_{1,1})=va^\equiv_2(K_{1,2})=1$ 
and $va^\equiv_2(K_{1,3})=va^\equiv_2(K_{2,2})=2.$  
If $m+n =4b+c$ where $b$ is a positive integer and $0\leq c \leq 3,$ 
then we have the following.\\ 
(1) For $c=0,$ if there are positive integers $h$ and $k$ such that 
$(m,n)=(4h,4k),$ then $va^\equiv_2(K_{m,n})=p(b: m,n),$ 
otherwise $va^\equiv_2(K_{m,n})= b+1.$\\  
(2) For $c=1,$ if there are positive integers $h$ and $k$ such that 
$(m,n)=(4h+1,4k)$ or $(4h,4k+1),$ then $va^\equiv_2(K_{m,n})=p(b: m,n),$ 
otherwise $va^\equiv_2(K_{m,n})= b+1.$\\ 
(3) For $c=2,$ if there are positive integers $h$ and $k$ such that 
$(m,n)=(4(h+1)+2, 4k),(4h+1,4k+1),$ or $(4h, 4(k+1)+2),$ 
then $va^\equiv_2(K_{m,n})=p(b: m,n),$ 
otherwise $va^\equiv_2(K_{m,n})= b+1.$\\ 
(4) For $c=3,$ if $(m,n)=(5,6)$ or 
there are positive integers $h$ and $k$ such that 
$(m,n)=(4(h+2)+3,4k),(4(h+1)+2, 4k+1),$ 
$(4h+1, 4(k+1)+2),$ or $(4h, 4(k+2)+3),$ 
then $va^\equiv_2(K_{m,n})=p(b: m,n),$ 
otherwise $va^\equiv_2(K_{m,n})= b+1.$ 
\end{Th}

\begin{Pf} 
It is easy to see that $va^\equiv_2(K_{1,1})=va^\equiv_2(K_{1,2})=1$ 
and $va^\equiv_2(K_{1,3})=va^\equiv_2(K_{2,2})=2.$  
Now consider the part $m+n =4b+c$ 
where $b$ is a positive integer and $0\leq c \leq 3.$ 
Since $(m+n)/b \geq 4,$ Lemma \ref{proper} yields that 
$K_{m,n}$ has an equitable $(b,2)$-tree-coloring if and only if  
$K_{m,n}$ has a proper equitable $b$-coloring. 
Thus each color class from an equitable $(b,2)$-tree-coloring 
of $K_{m,n}$ is an independent set. 

CASE 1: $c=0.$ An equitable $(b,2)$-tree-coloring of $K_{m,n}$ 
yields $b$ color classes of size 4. 
By Lemma $\ref{proper},$ each color class is independent. 
That is each color class is in a partite set $X$ or $Y.$  
This can happen if and only if there are 
positive integers $h$ and $k$ such that 
$(m,n)=(4h,4k).$ 

CASE 2: $c=1.$ An equitable $(b,2)$-tree-coloring of $K_{m,n}$ 
yields $b-1$ color classes of size 4 and 1 color class of size 5. 
By Lemma $\ref{proper},$ each color class is independent.  
That is each color class is in a partite set $X$ or $Y.$ 
This can happen if and only if there are 
positive integers $h$ and $k$ such that 
$(m,n)=(4h+1,4k)$ or $(4h,4k+1).$

CASE 3: $c=2.$ 

Subcase 3.1: $b=1.$ Then $m+n=6.$ 
One can see that $va^\equiv_2(K_{m,n})= 2.$ 

Subcase 3.2: $b \geq 2.$ 
An equitable $(b,2)$-tree-coloring of $K_{m,n}$ 
yields $b-2$ color classes of size 4 and 2 color classes of size 5.  
By Lemma $\ref{proper},$ each color class is independent. 
That is each color class is in a partite set. 
This can happen if and only if there are 
positive integers $h$ and $k$ such that 
$(m,n)=(4(h+1)+2, 4k),(4h+1,4k+1),$ or $(4h, 4(k+1)+2).$ 

CASE 4: $c=3.$ 

Subcase 4.1: $b=1.$ Then $m+n=7.$ 
One can see that $va^\equiv_2(K_{m,n})= 2=b.$ 

Subcase 4.2: $b=2.$ Then $m+n=11.$ 
Lemma \ref{bi0} yields that $K_{m,n}$ has 
an equitable $(q,2)$-tree-coloring for every $q \geq b+1=3.$ 
On the other hand, an equitable $(b,2)$-tree-coloring 
(that is an equitable $(2,2)$-tree-coloring) of $K_{m,n}$ 
yields 1 color class of size 5 and 1 color class of size 6.  
By Lemma $\ref{proper},$ each color class is independent. 
This can happen if and only if $(m,n)=(5,6).$ 
  
Subcase 4.3: $b\geq 3.$ 
An equitable $(b,2)$-tree-coloring of $K_{m,n}$ 
has $b-3$ color classes of size 4 and 3 color classes of size 5.  
By Lemma $\ref{proper},$ each color class is independent. 
This can happen if and only if there are 
positive integers $h$ and $k$ such that 
$(m,n)=(4(h+2)+3,4k),(4(h+1)+2, 4k+1),$ 
$(4h+1, 4(k+1)+2),$ or $(4h, 4(k+2)+3).$ 

Combining these facts with Lemma \ref{bi1}, 
we complete the proof.  
\end{Pf}

\section{$va^\equiv_2(K_{l,m,n})$} 

\begin{lem}\label{tri0123} 
Let $l+m+n =4b+c$ where $b$ is a positive integer. 
If $c \leq 2,$ then 
$K_{l,m,n}$ has an equitable $(t,2)$-tree-coloring for each $t\geq b+1.$   
If $c = 3,$ then 
$K_{l,m,n}$ has an equitable $(t,2)$-tree-coloring for each $t\geq b+2.$   
\end{lem} 
\begin{Pf}
For $c \leq 2,$ the proof is similar to that of Lemma \ref{bi0}. 
Now we assume $c =3.$ 
Let $m+n=4b+3 =rk+s(k+1)$ where $r$ is a positive integer,  
$s$ and $k$ are nonnegative integers. 
First consider the case $r+s=b+2.$ 
Then (i) $k\leq 2$  or (ii) $r=5, s=b-3, k=3$ and $b \geq 3.$ 
Again we can use Algorithm A to obtain an equitable $(b+1,2)$-tree-coloring. 

Finally, consider the case that $r+s \geq b+3.$ 
Then (i) $k\leq 2$  or (ii) $k=3$ and  $r \geq 5.$  
Again we can use Algorithm A  to 
obtain  an equitable $(r+s,2)$-tree-coloring. 
This completes the proof. 
\end{Pf}

\begin{lem}\label{trib012} 
Assume that $l+m+n =4b+c$ where $b$ is a positive  integer and $0\leq c \leq 2.$ 
If $K_{l,m,n}$ has an equitable $(b,2)$-tree-coloring,  
then $va^\equiv_2(K_{l,m,n})=p(b: l,m,n),$ 
otherwise $va^\equiv_2(K_{l,m,n})=b+1.$ 
\end{lem}
\begin{Pf} 
From Lemma \ref{tri0123}, 
$K_{l,m,n}$ has an equitable $(t,2)$-tree-coloring for each $t\geq b+1.$ 
By definition of $va^\equiv_2(K_{l,m,n})$, 
we have $K_{l,m,n}$ has no equitable $(b,2)$-tree-coloring 
if and only if  $va^\equiv_2(K_{l,m,n})=b+1.$ 

Assume $K_{l,m,n}$ has an equitable $(b,2)$-tree-coloring. 
Then each color class has size at least 4. 
By Lemma \ref{proper}, such equitable  $(b,2)$-tree-coloring 
is a proper equitable $b$-coloring. 
Thus $va^\equiv_2(K_{l,m,n})=p(b: l,m,n)$ by Lemma \ref{L3}.  
If $va^\equiv_2(K_{l,m,n})=p(b: l,m,n),$ 
then $K_{l,m,n}$ has an equitable $(b,2)$-tree-coloring by 
the definition of $va^\equiv_2(K_{l,m,n}).$ This completes the proof. 
\end{Pf}

\begin{Th}\label{triM012} 
If $l+m+n =4b+c$ where $b$ is a positive integer and $0\leq c \leq 2,$ 
then we have the following.\\
(1) For $c=0,$ if there are positive integers $j,h,$ and $k$ such that 
$(l,m,n)=(4j,4h,4k),$ then $va^\equiv_2(K_{l,m,n})=p(b: l,m,n),$ 
otherwise $va^\equiv_2(K_{l,m,n})= b+1.$\\  
(2) For $c=1,$ if there are positive integers $j,h,$ and $k$ such that 
$(l,m,n)=(4j+1,4h,4k), (4j,4h+1,4k),$ or $(4j,4h,4k+1),$ 
then $va^\equiv_2(K_{l,m,n})=p(b: l,m,n),$ 
otherwise $va^\equiv_2(K_{l,m,n})= b+1.$\\ 
(3) For $c=2,$ if there are positive integers $j,h,$ and $k$ such that 
$(l,m,n)= (4(j+1)+2, 4h, 4k),$ $(4j, 4(h+1)+2, 4k), (4j, 4h, 4(k+1)+2),$ 
$(4j+1, 4h+1, 4k),(4j+1, 4h,4k+1),$ or $(4j, 4h+1, 4k+1),$ 
then $va^\equiv_2(K_{l,m,n})=p(b: l,m,n),$ 
otherwise $va^\equiv_2(K_{l,m,n})= b+1.$
\end{Th}

\begin{Pf} 
Since $(l+m+n)/b \geq 4,$ Lemma \ref{proper} yields that 
$K_{l,m,n}$ has an equitable $(b,2)$-tree-coloring if and only if  
$K_{l,m,n}$ has a proper equitable $b$-coloring. 
Thus each color class from an equitable $(b,2)$-tree-coloring 
of $K_{l,m,n}$ is an independent set. 

CASE 1: $c=0.$ An equitable $(b,2)$-tree-coloring of $K_{l,m,n}$ 
yields $b$ color classes of size 4. 
By Lemma $\ref{proper},$ each color class is independent. 
That is each color class is in a partite set. 
This can happen if and only if there are 
positive integers $j, h,$ and $k$ such that 
$(l,m,n)=(4j,4h,4k).$ 

CASE 2: $c=1.$ An equitable $(b,2)$-tree-coloring of $K_{l,m,n}$ 
yields $b-1$ color classes of size 4 and 1 color class of size 5. 
By Lemma $\ref{proper},$ each color class is independent.  
That is each color class is in a partite set. 
This can happen if and only if there are 
positive integers $j, h,$ and $k$ such that 
$(l,m,n)=(4j+1,4h,4k), (4j,4h+1,4k),$ or $(4j,4h,4k+1).$  

CASE 3: $c=2.$ 

Subcase 3.1: $b=1.$ Then $l+m+n=6.$ 
One can see that $va^\equiv_2(K_{l,m,n})= 2.$ 

Subcase 3.2: $b \geq 2.$ 
An equitable $(b,2)$-tree-coloring of $K_{l,m,n}$ 
yields $b-2$ color classes of size 4 and 2 color classes of size 5.  
By Lemma $\ref{proper},$ each color class is independent. 
That is each color class is in a partite set. 
This can happen if and only if there are 
positive integers $j, h,$ and $k$ such that 
$(l,m,n)= (4(j+1)+2, 4h, 4k), (4j, 4(h+1)+2, 4k), (4j, 4h, 4(k+1)+2),$ 
$(4j+1, 4h+1, 4k),(4j+1, 4h,4k+1),$ or $(4j, 4h+1, 4k+1).$  

Combining these facts with Lemma \ref{trib012}, 
we complete the proof.  
\end{Pf}

\begin{Def}
We say that $(l,m,n)$ satisfies Condition A if 
there are positive integers $j,h,$ and $k$ 
such that $(l,m,n)=(4j,4h,4k -1), (4j ,4h -1,4k ),(4j -1,4h ,4k ),$ 
$(4j ,4h -2,4k -3), (4j ,4h -3,4k -2),  (4j -2,4h ,4k -3),$ 
$(4j -2,4h -3,4k ), (4j -3,4h ,4k -2),$ or  $(4j -3,4h -2,4k ).$ 
\end{Def}

\begin{lem}\label{tri3A} 
Let $l+m+n =4b+3$ where $b$ is a nonnegative integer. 
$K_{l,m,n}$ has an equitable $(b+1,2)$-tree-coloring 
if and only if $(l,m,n)$ satisfies condition A. 
\end{lem} 
\begin{Pf}
Assume that $G$ has an equitable $(b+1,2)$-tree-coloring. 
Then there are $b$ color classes of size 4 
and $1$ color class of size 3. 
By Lemma \ref{proper}, each color class of size 4 is independent. 
By definition of $(q,2)$-tree-coloring, 
a color class of size 3, say $C,$ is an independent set or 
$C$ induces $K_{1,2}.$ 

The case $C$ is an independent can happen 
if and only if there are positive integers $j,h,$ and $k$ 
such that$(l,m,n)=(4j,4h,4k-1), (4j,4h-1,4k),$ or $(4j-1,4h,4k).$ 

The case that $C$ induces $K_{1,2}$ can happen 
if and only if one element of $C$ is in one partite set 
and two other elements are in a different partite set. 
Thus the case that $C$ induces $K_{1,2}$ can happen 
if and only if there are positive integers $j,h,$ and $k$ 
such that $(4j,4h-2,4k-3), (4j,4h-3,4k-2),  (4j-2,4h,4k-3),$ 
$(4j-2,4h-3,4k), (4j-3,4h,4k-2),$ or  $(4j-3,4h-2,4k).$ 
This completes the proof. 
\end{Pf}

\begin{Def}
We say that $(l,m,n)$ satisfies Condition B if 
there are positive integers $j,h,$ and $k$ 
such that $(l,m,n)=(4(j+2)+3,4h,4k), (4j, 4(h+2)+3,4k), 
(4j,4h,4(k+2)+3,$ 
$(4(j+1)+2, 4h+1, 4k), (4(j+1)+2, 4h, 4k+1),$ 
$(4j+1,4(h+1)+2, 4k), (4j+1, 4h, 4(k+1)+2),$ 
$(4j, 4(h+1)+2, 4k+1),$ $(4j, 4h+1, 4(k+1)+2),$ or 
$(4j+1,4h+1,4k+1).$    
\end{Def} 

\begin{lem} \label{tri3B} 
Assume that $l+m+n =4b+3$ where $b$ is a positive integer. 
$K_{l,m,n}$ has an equitable $(b,2)$-tree coloring 
if and only if $(l,m,n)$ satisfies condition B. 
\end{lem} 

\begin{Pf}
Assume that $G$ has an equitable $(b,2)$-tree-coloring. 
Then there are $(b-3)$ color classes of size 4 
and $3$ color classes of size 5. 
By Lemma \ref{proper}, each color class  is independent. 
This can happen if and only if $(l,m,n)$ satisfies condition B. 
\end{Pf}

\begin{lem}\label{trib3} 
Let $l+m+n =4b+3$ where $b$ is a positive  integer. 
We have the following.\\ 
(1) $K_{l,m,n}$ has no equitable $(b+1,2)$-tree-coloring 
if and only if  $va^\equiv_2(K_{l,m,n})=b+2.$\\ 
(2) Assume that $K_{l,m,n}$ has an equitable $(b+1,2)$-tree-coloring. 
If $K_{l,m,n}$ has an equitable $(b,2)$-tree-coloring,   
then  $va^\equiv_2(K_{l,m,n})=p(b: l,m,n),$ 
otherwise $va^\equiv_2(K_{l,m,n})=b+1.$  
\end{lem}

\begin{Pf}
From Lemma \ref{tri0123}, 
$K_{l,m,n}$ has an equitable $(t,2)$-tree-coloring for each $t\geq b+2.$ 
By definition of $va^\equiv_2(K_{l,m,n})$, 
we have $K_{l,m,n}$ has no equitable $(b+1,2)$-tree-coloring 
if and only if  $va^\equiv_2(K_{l,m,n})=b+2.$ 

Now assume that $K_{l,m,n}$ has an equitable $(b+1,2)$-tree-coloring. 
Thus $K_{l,m,n}$ has an equitable $(t,2)$-tree-coloring for each $t\geq b+1.$ 
If $K_{l,m,n}$ has no equitable $(b,2)$-tree-coloring, 
then   $va^\equiv_2(K_{l,m,n})=b+1$ by the definition. 
Consider the case that $K_{l,m,n}$ has an equitable $(b,2)$-tree-coloring. 
Thus each color class has size at least 4. 
By Lemma \ref{proper}, such a coloring is an equitable $b$-coloring. 
Lemma \ref{L3} yields $va^\equiv_2(K_{l,m,n})=p(b: l,m,n).$ 
\end{Pf}

\begin{Th}\label{triM3} 
$va^\equiv_2(K_{1,1,1})=2.$  
Assume that $l+m+n =4b+3$ where $b$ is a positive integer. 
Then we have the following.\\
(i) If $(l,m,n)$ does not satisfy Condition A, 
then $va^\equiv_2(K_{l,m,n})=b+2.$\\  
(ii) If $(l,m,n)$ satisfies Condition A but does not satisfy Condition B, 
then $va^\equiv_2(K_{l,m,n})=b+1.$\\  
(iii) If $(l,m,n)$ satisfies Condition A and Condition B, 
then $va^\equiv_2(K_{l,m,n})=p(b+1:l,m,n).$\\  
\end{Th} 
\begin{Pf} 
It is easy to see that $va^\equiv_2(K_{1,1,1})=2.$  
Now consider the part $l+m+n =4b+3$ with a positive integer $b.$ 
Using Lemmas \ref{tri3A} and \ref{trib3} (1), we have (i). 
Using Lemmas \ref{tri3A}, \ref{tri3B}, and \ref{trib3} (2), we have (ii) and (iii).  
This completes the proof. 
\end{Pf}



\begin{thebibliography}{99}
\bibitem{CLW94}
B.-L. Chen, K.-W. Lih, and P.-L. Wu, Equitable coloring and the
maximum degree, \emph{ Europ. J. Combinatorics} 15(1994), 443--447. 

\bibitem{CY12}
B. L. Chen, C. H. Yen, 
Equitable $\Delta$-coloring of graphs, 
\emph{ Discrete Math.} 312(2012) 1512--1517.

\bibitem{Fan11} 
H. Fan, H. A. Kierstead, G. Z. Liu, T. Molla, J. L. Wu, and X. Zhang, 
A note on relaxed equitable coloring of graphs, 
\emph{ Inform. Process. Lett.} 111(2011) 1062--1066.

\bibitem{GareyJohnson}
M. R. Garey and D. S. Johnson, \emph{ Computers and Intractability: A
Guide to the Theory of NP-completeness,} W. H. Freeman and Company,
New York, 1979.

\bibitem{HS}
A. Hajnal and E. Szemer\' edi,
Proof of conjecture of Erd\H os,
in: \emph{ Combinatorial Theory and its Applications, Vol. II}
(P.~Erd\H os, A.~R\' enyi and V.~T.~S\' os Editors),
(North-Holland, 1970), 601--623.

\bibitem{KK08}
H. A. Kierstead and A. V. Kostochka,
A short proof of the Hajnal-Szemer\' edi Theorem on equitable colouring,
\emph{ Combin. Probab. Comput.} 17(2008), 265--270.

\bibitem{Ko}
A. V. Kostochka,
Equitable colorings of outerplanar graphs,
\emph{ Discrete Math.} 258(2002), 373--377.

\bibitem{LiBu09}
Q. Li and Y. Bu,
Equitable list coloring of planar graphs without 4- and 6-cycles,
\emph{ Discrete Math.} 309(2009),  280--287.

\bibitem{LW}
K.-W. Lih and P.-L. Wu,
On equitable coloring of bipartite graphs,
\emph{ Discrete Math.} 151(1996),  155--160.

\bibitem{M}
W. Meyer,
Equitable Coloring,
\emph{ American Math. Monthly}
 80(1973), 920--922.

\bibitem{Nak12}
K. Nakprasit,
Equitable  colorings of planar graphs with maximum degree at least nine,
\emph{ Discrete Math.} 312(2012), 1019--1024. 

\bibitem{NN12}
K. Nakprasit and K. Nakprasit, 
Equitable colorings of planar graphs without short cycles, 
\emph{ Theoretical Computer Science} 465 (2012) 21--27.

\bibitem{NN1ab}
K. M. Nakprasit and K. Nakprasit, 
The strong equitable vertex 1-arboricity of complete multipartite graphs, 
\emph{ Manuscript}  (2015).

\bibitem{west}
D. B. West,  Introduction to graph theory, Second edn, Printice Hall, pp. xx+588, 2001.

\bibitem{WZL13} 
J.-L. Wu, X. Zhang, and H.L. Li, 
Equitable vertex arboricity of graphs, \emph{ Discrete Math.}
313(2013), 2696--2701.

\bibitem{YW12}
Z. Yan and W. Wang, 
Equitable chromatic threshold of complete multipartite graphs, 
arXiv:1207.3578 [math.CO]. 


\bibitem{YZ1}
H.-P. Yap and Y. Zhang,
The equitable $\Delta$-colouring
conjecture holds for outerplanar graphs,
\emph{ Bull. Inst. Math. Acad. Sin.}
25(1997), 143--149. 

\bibitem{ZW11} 
X. Zhang, J.-L. Wu, On equitable and equitable list colorings of series-parallel graphs, 
\emph{ Discrete Math.} 311 (2011) 800–803.

\bibitem{ZY98}
Y. Zhang and H.-P. Yap,
Equitable colourings of planar graphs,
\emph{ J. Combin. Math. Combin. Comput.}
27(1998),  97--105.

\bibitem{ZhuBu08}
J. Zhu and Y. Bu,
Equitable list colorings of planar graphs without short cycles,
\emph{ Theoretical Computer Science} 407(2008), 21--28.

\end{thebibliography}
\end{document}